# A fast semi-analytic algorithm for computing solutions associated with multiple moving or fixed bottlenecks: Application to joint scheduling and signal timing

Michele D. Simoni[a], Christian G. Claudel[b],*

[a]Dept. of Civil, Architectural, and Environmental Engineering, University of Texas at Austin, 301E Dean Keeton, St., Austin, 78712 TX, USA
Dept. of Civil, Architectural, and Environmental Engineering, University of Texas at Austin, 301E Dean Keeton, St., Austin, 78712 TX, USA

**Abstract**

Moving and fixed bottlenecks are moving or fixed capacity restrictions that affect the propagation of traffic flow. They are a very important modeling approach to describe the effects of slow vehicles and traffic signals in transportation networks. However, the computation of solutions associated with the presence of fixed and moving bottlenecks is complex, since they both influence and are influenced by traffic. In this study, we propose a fast numerical scheme that can efficiently compute the solutions to an arbitrary number of fixed and moving bottlenecks, for a stretch of road modeled by the Lighthill-Whitham-Richards (LWR) model with triangular fundamental diagram. The numerical scheme is based on a semi-analytic Lax-Hopf formula that requires a very low number of operations compared with existing schemes. We illustrate the performance of the numerical scheme on scenarios involving multiple slow vehicles and traffic signals, and demonstrate that this scheme can be part of an optimization loop to simultaneously optimize the schedule of several heavy-duty vehicles and traffic signals in a city for alleviating traffic congestion

## 1. Introduction

In traffic flow theory, different typologies of "slow" vehicles (or platoons) can be modeled as "moving bottlenecks". These obstructions in traffic streams are usually associated with the presence of buses in urban traffic, and trucks or simply slower vehicles on highways. All these situations, indeed, are characterized by a partial blockage the road (typically the right lane in right hand driving countries), causing a capacity reduction. The concept of moving bottleneck can be extended to fixed bottlenecks, which represent static (spatially) and time varying capacity restrictions caused for example by traffic lights and traffic incidents.

Some of the main challenges of modeling moving bottlenecks consist of identifying and modeling features regarding their speed (depending on the traffic conditions and on the maximum speed of the vehicle), their discharging flow (maximum rate at which vehicles overtake) and the entity of queue held back. Several studies have highlighted the importance of the effects of moving bottlenecks on traffic (Munoz and Daganzo, 2002; Daganzo and Laval, 2005) and have developed methodologies to include them into existing traffic models. Gazis and Herman developed in 1992 a model based on the conservation of flow, unconditional existence of the flow-density relation, and independence of capacity state from the bottleneck state. The first complete formulation based on the Lighthill–Whitham–Richards (LWR) model was proposed few years later by Newell (1993; 1998) where the moving bottleneck is assumed to behave as in a scaled-down version of the freeway's fundamental diagram not influenced by the bottleneck speed. In recent years, more comprehensive formulations of the moving bottleneck problem have been proposed by Munoz and Daganzo (2002), Leclerq et al. (2004) and Daganzo and Laval, (2005). Other studies have focused on numerical methods to solve the fixed and moving bottleneck problems within the LWR model (Lebacque et al., 1998; Giorgi et al., 2002; Leclercq, 2007). Kerner and Klenov (2010) explored thoroughly features of moving bottlenecks such as the critical speed at which traffic breakdown, based on "Three-phase traffic theory".

* Corresponding author. Tel.: +1 512 705 7195;
  *E-mail address:* christian.claudel@utexas.edu

In recent years, moving bottlenecks were also studied in the field of applied mathematics by Lattanzio et al. (2011) and Gasser et al. (2013), who developed a coupled PDE-ODE model to reproduce the dynamics between traffic flows and moving bottlenecks like buses.

In this article, we similarly propose an approach that accounts for both the impacts of moving bottlenecks on surrounding traffic and the converse. In addition, we propose an efficient algorithm allowing the simulation of an arbitrary number of moving bottlenecks associated with different maximum speeds.

To achieve this, we propose in this article a new efficient formulation that computes the parameters associated with moving and fixed bottlenecks (trajectories and passing flows), without having to compute the complete solution, thereby improving computational times by orders of magnitude over classical numerical schemes, without affecting the computational accuracy.

As mentioned earlier, the problem of computing the trajectories and parameters (passing flows) associated with moving bottlenecks is not straightforward, since the bottlenecks both influence and are influenced by surrounding traffic. Thus, in order to compute the density map associated with a general problem (involving initial, boundary conditions and bottlenecks), it is necessary to simultaneously compute the solution to the LWR model and the corresponding trajectories of the bottlenecks, since these are initially unknown. This process is computationally intensive. In particular, this requires us to map the solution on the entire computational domain, since the trajectories of the moving bottlenecks are affected by the solution itself.

The algorithm we propose instead allows us to determine the parameters and trajectories of the moving bottlenecks without requiring us to determine the solution on the entire computational domain. It is based on an extension of the semi-analytical solutions to arbitrary Hamilton-Jacobi equations introduced in (Mazaré et al., 2011). Using semi-explicit solutions, we show that the trajectories of an arbitrary number of fixed and moving bottlenecks can be simultaneously marched forward in time for a very low computational cost. Indeed, if the piecewise affine initial conditions contain $n_i$ blocks, the piecewise affine upstream and downstream boundary conditions contain $n_u$ and $n_d$ blocks respectively, and $n_b$ bottlenecks are considered, the future evolution of each bottleneck can be computed by at most ($n_u+n_b+2$) calculations of explicit functions, which determine the future value of the solution to the Hamilton-Jacobi equation along the trajectory of the bottleneck. Once this set of calculations is done, the future evolution of the moving bottleneck is completely determined, in function of the difference between the current value of the solution to the Hamilton Jacobi equation along the trajectory, and its future value along the predicted trajectory. This process is marched forward in time, and allows one to simultaneously compute the parameters associated with all moving and fixed bottlenecks of the problem, without having to compute the solution everywhere (indeed, solutions are only required along the trajectories of the bottleneck, greatly reducing the computational time required to solve the problem).

Once the parameters and trajectories of all moving and fixed bottlenecks are known, one can use this information to efficiently compute the solution of the problem everywhere using the Lax-Hopf algorithm (which was shown to be faster than the Godunov scheme if solutions are only required at a given time horizon in (Claudel and Bayen, 2010a). Since the Lax-Hopf algorithm can compute the solution at any point of the space time domain using only initial, boundary and bottleneck data, this approach is well adapted to optimization problems in which we are only interested in knowing the solutions at a limited number of points (on which the objective function of the problem depend).

Another advantage of this algorithm is its very favorable computational error characteristics. The only errors induced by the proposed scheme are errors related to the discretization in time of the fixed and moving bottleneck trajectories and passing flows, and an approximation of the behavior of the bottlenecks around intersections of bottleneck trajectories (if such intersections occur). No other numerical error is induced since the solutions to initial, boundary and bottleneck condition blocks are explicit and exact. Since all other non-event-based numerical methods to solve moving bottleneck problems (for example based on LWR (Leclercq, 2007) or on the Variational Method (Daganzo and Laval, 2005) also require discretized moving and fixed bottleneck trajectories but also use approximate solution methods to solve the LWR equation, the proposed algorithm yields more accurate solutions than such methods. Event-based methods can be exact, though they still require the computation of the solution on the entire computational domain, and to date no wave-front tracking formulation capable of handling multiple moving bottlenecks exist (Delle Monache and Goatin, 2014).

Thanks to these favorable properties, the proposed algorithm could be used to efficiently tackle complex traffic estimation and control problems characterized by presence of multiple trucks or buses. For this reason, we present in this study an application of this framework to the optimization of multiple traffic signals times and buses schedules, solved by coupling the present algorithm to a genetic algorithm-based optimization framework.

In the remainder of this article, we first introduce the background theory adopted in this study for the modeling of moving bottlenecks. Then, we provide a description of the fast semi-analytic algorithm to simulate single moving bottlenecks. We then extend this algorithm to scenarios including multiple moving and fixed bottlenecks, which can be associated with different types of vehicles. Finally, we illustrate the application of the proposed algorithm to different optimization problems. We conclude with a number of general remarks and recommendations for future research.

## 2. Analytical solutions to the Hamilton-Jacobi PDE

In this Section we briefly summarize the main features of the macroscopic traffic simulation used to investigate moving bottlenecks. The LWR model and the Hamilton-Jacobi PDE are described respectively in Section 2.1 and Section 2.2. The generalized Lax-Hopf formula used to solve this problem is presented in Section 2.3 and the formulation of initial, boundary and internal conditions is provided in Section 2.4. Finally, in Section 2.5 we describe the model for moving bottlenecks adopted for the derivation of the internal conditions.

### 2.1. The LWR-PDE

Given a one-dimensional uniform section of highway, limited by $x_0$ upstream and $x_n$ downstream. For a given time $t$ and position $x$ we define the local traffic density $k(x,t)$ in vehicles per unit of length, and the instantaneous flow $q(x,t)$ in vehicles per unit time. The conservation of vehicles on the highway is written as follows (Lighthill and Whitman, 1956; Richards, 1956; Garavello and Piccoli, 2006).

$$\frac{\partial k(t,x)}{\partial t} + \frac{\partial q(t,x)}{\partial x} = 0 \qquad (1)$$

For first order traffic flow models, flow and density are related by the Fundamental Diagram (FD); in this article we adopt triangular FD (Daganzo, 1994). The FD is a positive function defined on $[0,kj]$, where $k_j$ is the maximal density (jam density). It ranges in $[0,q_{max}]$ where $q_{max}$ is the maximum flow (capacity). It is assumed to be differentiable with derivative $Q'(0)=v>0$ (free flow speed) and $Q'(k\_j)=w<0$ (congested wave speed), and it is defined as follows:

$$Q(k) = \begin{cases} v\,k & : \ 0 \leq k \leq k_c \\ -w\,(k - k_c) & : \ k_c \leq k \leq k_j \end{cases} \qquad (2)$$

*2.2. The Moskovitz function*

The Moskovitz function expresses the cumulated vehicle count *N(x,t)* and it represents the continuous vehicle count at location *x* and time *t*. In the Moskovitz framework one assumes that all vehicles are labeled by increasing integers as they pass the entry point $x_0$ of a highway section, and that they cannot pass each other. If the latest car that passed an observer standing at location *x* and time *t* is labeled n, then *N(x,t)=n*.

Replacing *k* and *q* with *N* yields to Hamilton-Jacobi PDE (Newell, 1993; Daganzo, 2005a, 2006; Claudel and Bayen, 2010a):

$$\frac{\partial N(x,t)}{\partial t} - Q\left(-\frac{\partial N(x,t)}{\partial x}\right) = 0 \tag{3}$$

*2.3. The generalized Lax-Hopf Formula*

From Aubin et al. (2008), the solution associated with the value condition function c, denoted by N_c, is the infimum of an infinite number of functions of the value condition:

$$N_c = \inf\{c(t-T, x-Tu) + TR(u)\} \; s.t. \; (u,T) \in [w, v_f] \times R_+ \; and \; (t-T, x-Tu) \in Dom(c) \tag{4}$$

Where c(x,t) corresponds to:

$$c(x,t) = \begin{cases} N_{ini}(x) & t = 0 \\ N_{up}(t) & x = x_0 \\ N_{down}(t) & x = x_n \end{cases} \tag{5}$$

And $R(u)$, which is convex transform associated with the fundamental diagram:

$$R(u) = \sup_{k \in [0, k_j]} (Q(k) - u \cdot k) \tag{6}$$

This equation is well known in the Hamilton-Jacobi literature and often referred to as Lax-Hopf formula (Aubin et al., 2008; Evans, 1998).
Assuming a triangular fundamental diagram, the calculation of its convex transform R yields to:

$$\forall u \in [w, v_f], R(u) = k_c(v_f - u) \tag{7}$$

*2.4. Boundary and internal conditions based on triangular fundamental diagram*

a) <u>Definition of initial, upstream, downstream and internal conditions</u>
The initial condition can be expressed as a piecewise linear function, with each linear piece defined by:

$$c_{ini_i}(x) = \begin{cases} -k_i x + b_i & : x_i \leq x \leq x_{i+1} \\ +\infty & : otherwise \end{cases} \tag{8}$$

With the above definition, the initial condition can be written as $c_{ini} = \min_i c_{ini_i}$

Similarly, the upstream boundary condition is assumed to be piecewise linear, with each piece defined by:

$$c_{up_j}(t) = \begin{cases} q_j t + d_j & : t_j \leq t \leq t_{j+1} \\ +\infty & : otherwise \end{cases} \quad (9)$$

With this definition, the upstream boundary condition can be written as $c_{up} = \min_j c_{up_j}$

The downstream boundary condition is also assumed to be a piecewise linear function, with each piece defined by:

$$c_{down_j}(t) = \begin{cases} p_j t + b_j & : t_j \leq t \leq t_{j+1} \\ +\infty & : otherwise \end{cases} \quad (10)$$

This enables us to define the downstream boundary condition function as $c_{down} = \min_j c_{down_j}$,

The internal condition corresponding to a fixed or moving bottleneck active on the domain can be defined as:

$$c_{int}(t,x) = \begin{cases} N_b + \dfrac{(N_e - N_b)}{(t_e - t_b)} \cdot (t - t_b) & : x = x_b + \dfrac{(x_e - x_b)}{t_e - t_b} \cdot (t - t_b) \text{ and } t \in [t_b, t_e] \\ +\infty & : \qquad\qquad otherwise \end{cases} \quad (11)$$

One of the major results of Mazaré et al. (2011) is that the solutions associated with each linear piece of the initial, upstream, downstream and internal boundary conditions can be computed analytically as follows:

b)  Solution to a linear initial condition
If $0 \leq k_i \leq k_c$, the initial condition imposes a free-flow state.

$$N_{c_{ini}}(x,t) = \begin{cases} k_i(tv_f - x) + b_i & : x_i + tv_f \leq x \leq x_{i+1} + tv_f \\ k_i(tv_f - x) + b_i + x_i(k_c - k_i) & : x_i + tw \leq x \leq x_{i+1} + tv_f \end{cases} \quad (12)$$

else, if $k_c \leq k_i \leq k_j$:

$$N_{c_{ini}}{}^i(x,t) = \begin{cases} k_i(tw - x) - t k_j w + b_i & : x_i + tw \leq x \leq x_{i+1} + tw \\ k_c(tw - x) - t k_j w + x_{i+1}(k_c - k_i) + b_i & : x_{i+1} + tw \leq x \leq x_{i+1} + tv_f \end{cases} \quad (13)$$

c)  Solution to a linear upstream boundary condition
For an upstream boundary condition $N_{up}$ defined as: $N_{up}{}^j(t) = q_j t + d_j$ with $d_j = -q_j t + \sum_{l=0}^{j-1}(t_{l+1} - t_l) q_j{}^l$, the solution component can be expressed as:

$$N_{c_{up}{}^j}(x,t) = \begin{cases} d_j + q_j\left(t - \dfrac{x - x_0}{v_f}\right) & : x_0 + v_f(t - t_{j+1}) \leq x \leq x_0 + v_f(t - t_j) \\ d_j + q_j t_{j+1} + k_c\left((t - t_{j+1})v_f - (x - x_0)\right) & : x_0 \leq x \leq x_0 + v_f(t - t_{j+1}) \end{cases} \quad (14)$$

d)  Solution to a linear downstream boundary condition
For a downstream boundary condition $N_{down}{}^j$, defined as $N_{down}{}^j(t) = p_j t + b_j$ with $b_j = -p_j t + N_{ini}^{(n-1)}(x_n) + \sum_{l=0}^{j-1}(t_{l+1} - t_l) q_j{}^l$, the solution component can be expressed as:

$$N_{down}{}^j(x,t) = \begin{cases} b_j + p_j t - \left(\frac{p_j}{w} + k_j\right)(x_n - x) & : x_n + w(t - t_j) \leq x \leq x_n + w(t - t_{j+1}) \\ b_j + p_j t_{j+1} + k_c\left((t - t_{j+1})v_f + x_n - x\right) & : x_n + w(t - t_j) \leq x \leq x_n \end{cases} \quad (15)$$

e) <u>Solution to a linear internal condition</u>

For an internal condition $N_{int}$, the solution component can be expressed as:

$$N_{int}(x,t) = (t - t') \cdot (u + v) \cdot k_c + (N_e - N_b) \cdot \frac{(t' - t_b)}{(t_e - t_b)} + N_b \quad : x \leq x_b + v \cdot (t - t_b) \wedge x \geq x_b + w \cdot (t - t_b) \wedge t \geq t_b \quad (16)$$

where $t'$, which corresponds to the capture time in the viability framework from which these formulations are derived (Aubin, et al., 2008), is derived as follows:

$$t' = \begin{cases} t - \frac{(x_b + s \cdot (t - t_b) - x)}{(s - v)} & : x_e + v \cdot (t - t_e) \leq x \wedge x_b + s \cdot (t - t_b) \leq x \\ t_e & : x_e + w \cdot (t - t_e) \leq x \\ t - \frac{(x_b + s \cdot (t - t_b) - x)}{(s + w)} & : else \end{cases} \quad (17)$$

and $u$ corresponds to associated optimal control of the auxiliary dynamical system (Aubin et al., 2008; Bayen et al., 2007):

$$u = \begin{cases} -v & : x_e + v \cdot (t - t_e) \leq x \wedge x_b + s \cdot (t - t_b) \leq x \\ \frac{(x_e - x)}{(t - t_e)} & : x_e + w \cdot (t - t_e) \leq x \\ w & : else \end{cases} \quad (18)$$

Finally, $s$ corresponds to the velocity of the moving bottleneck:

$$s = \frac{(x_e - x_b)}{(t_e - t_b)} \quad (19)$$

and r corresponds to the number of vehicles passing the moving bottleneck per unit time:

$$r = \frac{N_e - N_b}{t_e - t_b} \quad (20)$$

These two parameters are unknown at the beginning of our problem. Indeed, only the initial position and the starting time of each moving bottleneck are known a priori, but not the evolution of the parameters s and r associated with each moving bottleneck, since they depend on the solution itself. As the matter of fact, the objective of the present article is to compute the evolution of s and r for each bottleneck, given known initial and boundary conditions, and given the knowledge of maximal velocity, starting time and starting position of the bottlenecks.

*2.5. Modeling single moving bottlenecks as internal conditions*

The dynamics of s and r is complex, since the behavior of a moving bottleneck is inherently hybrid, with active and inactive phases depending on the state of traffic. A moving bottleneck becomes "active" when it actually slows down the incoming traffic from upstream. This situation occurs when the traffic flow is sufficiently high to be hindered by the moving bottleneck. Following Munoz and Daganzo (2004), and Daganzo and Laval (2005), three situations can be distinguished (Figure 1):

1. The moving bottleneck is inactive because there is enough capacity for regular traffic to overtake (region 1 in Figure 1)

2. The moving bottleneck is active because regular traffic is traveling at higher speed and capacity is not enough for everyone to overtake (region 2 in Figure 1)

3. The moving bottleneck is inactive because regular traffic is traveling at a lower speed than the maximum velocity of the bottleneck, because of congestion (region 3 in Figure 1)

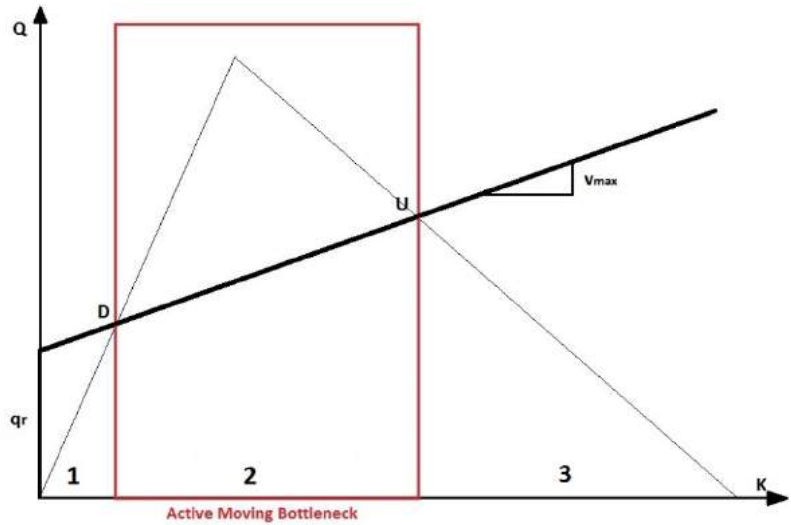

Figure 1: Flow-density relationship of moving bottlenecks according to the Munoz-Daganzo model

In order to identify whether the moving bottleneck is active and derive its corresponding internal conditions we adopted the following approach based on the difference of cumulated flow between two consecutive points along the trajectory of the moving bottleneck:

- Choose an arbitrary time step $\Delta t$
- Calculate the values of the Moskovitz function for: $M(x_0, t_0) = M_0$ and $M(x_0 + v_{max}\Delta t, t_0 + \Delta t) = M_1$, where $(x_0, t_0)$ corresponds to the position of the moving bottleneck in the end of the previous time interval, and $v_{max}$ corresponds to the maximum speed of the moving bottleneck.
- Identify the three abovementioned cases based on the flow between the two consecutive points (ratio between the difference of the Moskovitz function and time):

   a) $0 < \frac{(M_1 - M_0)}{\Delta t} < q_r$ → inactivity due to low volumes (traffic is too light)

   b) $\frac{(M_1 - M_0)}{\Delta t} > q_r$ → activity

   c) $\frac{(M_1 - M_0)}{\Delta t} < 0$ → inactivity due to congestion (traffic is slower than the maximum velocity of the bottleneck)

In the above, $q_r$ corresponds to the maximum passing rate of the moving bottleneck, which is the maximum flow that can ever pass the moving bottleneck going at its maximum speed. The formulation of $q_r$ is based on the model by Munoz and Daganzo (2004) and it corresponds to:

$$q_r = \frac{(u - v_{max}) \cdot k_c \cdot (n_l - 1)}{n_l} \tag{21}$$

where $u$ stands for the free flow speed, $k_c$ is the critical density and $n_l$ is the number of lanes.

Only in the case of active moving bottlenecks, a new internal condition with speed $s = v_{max}$ and overtaking rate $r = q_r$ is defined and stored. The internal condition applies between times $t_0$ and $t_0 + \Delta t$, and between positions $x_0$ and $x_0 + v_{max} \cdot \Delta t$, with beginning and end values of $M_0$ and $M_1 = M_0 + q_r \Delta t$.

In case of activity over several consecutive intervals ($\Delta t$), only the values of the Moskovitz function at the onset and end of activity, as well as the corresponding times and positions, are stored as internal conditions $c_{int}$.

In case of inactivity of the moving bottleneck due to congested conditions, the moving bottleneck travels at the speed of the surrounding traffic (which is less than its maximal speed), given by:

$$v = \frac{-w \cdot (k_0 - k_j)}{k_0} \tag{22}$$

Where $w$ corresponds to the congested speed, $k_j$ corresponds to the jam density and $k_0$ corresponds to the density of the traffic around the bottleneck.

- Finally, if the moving bottleneck is inactive due to low flow conditions, its velocity is set to $v_{max}$.

This algorithm can be summarized as the below pseudocode (Algorithm 1).

**Algorithm 1: Pseudo-code for the computation of internal conditions associated with a single moving bottleneck**

| | |
|---|---|
| **Input: $(x_1, t_1, x_2, v_{mb})$** moving bottleneck; | *Input initial position and time, final position, and performance characteristics of the moving bottleneck* |
| **Input: $(v, k_c, k_j, n_l)$** | *Input list of fundamental diagram parameters* |
| **Input T** | *Input simulation time horizon* |
| $q_r = (v - v_{mb}) * k_c * (n_l - 1)/n_l$; | *Derive maximum passing rate for the moving bottleneck, given the number of lanes of the road* |
| $t_0 = t_1$ <br> $x_0 = x_1$ | *Initialize bottleneck time* <br> *Initialize bottleneck position* |
| **while** $t_0 \leq T$ **and** $x_0$ | *While the bottleneck vehicle is still on the computational domain* |
| $t_1 = t_0 + \Delta t$ <br> $x_1 = x_0 + v_{mb}\Delta t$ | *update time* <br> *update position* |
| $M_0 = \inf\{M_{ini}, M_{up}, M_{down}, M_{int}\}$ calculated at point $(x_0, t_0)$ using initial, upstream, downstream and currently defined internal conditions <br><br> $M_1 = \inf\{M_{ini}, M_{up}, M_{down}, M_{int}\}$ calculated at point $(x_1, t_1)$ using initial, | *Calculate Moskovitz function at the previous and new positions of the moving bottleneck* |

*upstream, downstream and currently defined internal conditions*

---

**if $(M_1 - M_0)/\Delta t > 0$ then**
    **if $(M_1 - M_0)/\Delta t < q_r$**                       *If bottleneck is inactive due to low flows*

$t_0 = t_1$                                                                               *Update time*
$x_0 = x_1$                                                                              *Update position*

    **Else**                                                                                   *If bottleneck is active*

        **Add new internal condition with parameters**        *Store new internal condition*
                **$\{N_b, x_b, t_b; N_e, x_e, t_e\}$**

$t_0 = t_1$                                                                                 *Update time*
$x_0 = x_1$                                                                           *Update position*

    **end if**                                                                      *If bottleneck is inactive due to high congestion*

**Else**                                                                                *Compute actual speed of bottleneck*

**derive speed s from Equation $\frac{-w \cdot (k_0 - k_j)}{k_0}$**        *Update time*
                                                                                     *Update position*

$t_1 = t_0 + \Delta t$
$x_1 = x_0 + s\Delta t$

**end if**

**end while**

---

To illustrate the capabilities of Algorithm 1, we present the following example of a stretch of a two-lane road of length 3000 m characterized by some arbitrary initial and boundary shown respectively in Table 1 and Table 2. A moving bottleneck entering at $x = 1500\ m$ and $t = 150\ s$ with $v_{max} = 5\ m/s$ is included in the simulation (red trajectory in Figure 2). The following parameters characterizing the triangular fundamental diagram are chosen: $u = 30\ m/s$, $k_c = 0.04\ veh/m$, $k_j = 0.2\ veh/m$.

Table 1: upstream and downstream boundary conditions

| $i$ | $t_{i-1}$ | $t_i$ | $q_{up}^i$ | $q_{down}^i$ |
|---|---|---|---|---|
| 1 | 0 | 40 | 1.0 | 0.9 |
| 2 | 40 | 180 | 1.0 | 0.2 |
| 3 | 180 | 300 | 1.0 | 0.9 |

Table 2: initial conditions

| $i$ | $x_{i-1}$ | $x_i$ | $k_{ini}$ |
|---|---|---|---|
| 1 | 0 | 1000 | 0.04 |
| 2 | 1000 | 2000 | 0.02 |
| 3 | 2000 | 3000 | 0.04 |

The simulation, which was performed on Matlab with a 2.3 GHz processor, required about 0.75 seconds, of which less than 0.03 seconds are spent computing the parameters of the internal conditions and the trajectory of the moving bottleneck. The rest is used to compute the solution on the entire on a rectangular grid of resolution one second and ten meters. The results of the simulation are shown in a space-time-density diagram (Figure 2).

This simulation illustrates the benefits of the method over existing numerical schemes. Different numerical schemes have been proposed to model moving bottlenecks (such as first order numerical schemes (Daganzo and Laval, 2005; Leclercq, 2007), variational schemes (Daganzo, 2005), or wave-front tracking schemes (Henn, 2005), although they all require the solution to be computed everywhere on the computational domain. For example, the Godunov scheme (first order) requires us to compute the solution on the entire computational grid, and so does the Variational theory (in which bottlenecks are encoded as shortcuts in the computational grid). Similarly, wave-front tracking methods require the solution to be computed on the entire computational domain. Since most optimization problems in transportation do not require us to know the solution everywhere on the computational domain, this is a significant advantage as it allows us to first compute the parameters of all moving bottlenecks, and then compute the exact solution at the few points of the computational domain needed to determine the objective to be optimized, corresponding to a significant improvement in computational time and complexity.

## 3. Fast algorithm to compute multiple bottlenecks

In this section, after a brief review the moving bot we objective is to expand the previous algorithm to the case of an arbitrary number of moving and fixed bottlenecks, which can have distinct parameters (maximum velocity).

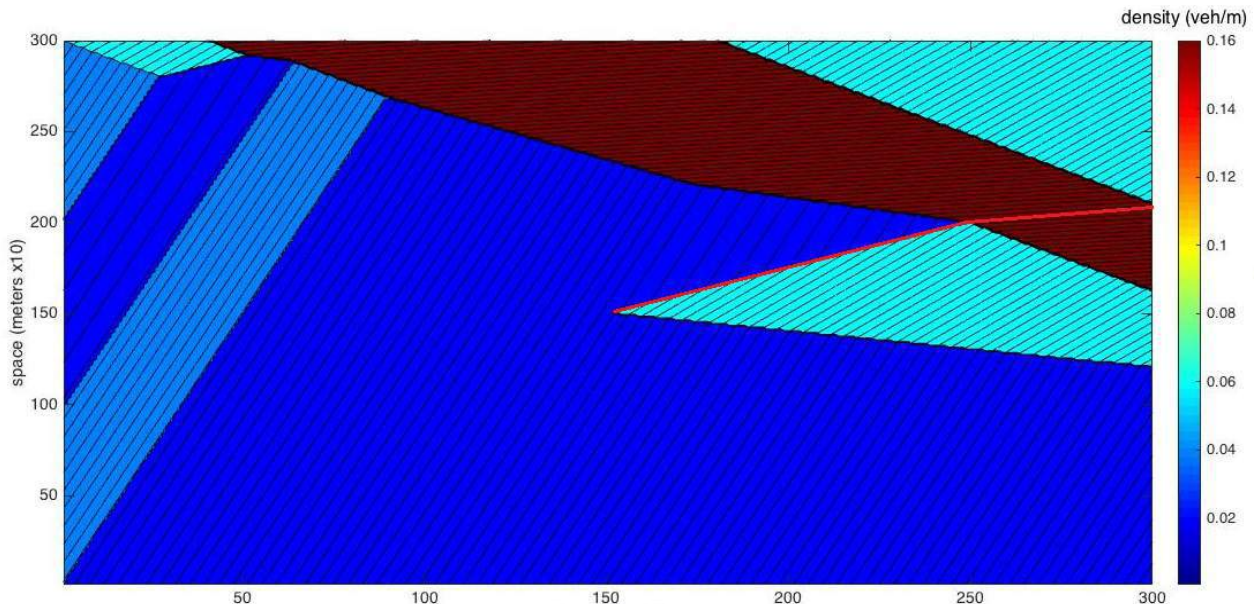

Figure 2: space-time-density diagram representing the results of the test stimulation. In this simulation, the trajectory of the moving bottleneck is shown in red

*3.1. Background*

The necessity of modeling the impacts of multiple trucks, buses and other kinds of slow moving vehicles on traffic has been recognized and increasingly emphasized in the last twenty years in the field of Traffic Flow Theory. To our knowledge, all these efforts have been made to study the effects of a single moving bottleneck, or moving bottlenecks for which the parameters (activity, velocity) are determined before the simulation. The objective is t include these into current traffic models, such as the LWR model. Extending previous work to an arbitrary number of moving bottlenecks (without relying on continuous flow approximations as in Liu, et al., 2015) implies several challenges, consisting in dynamically computing several trajectories-since moving bottlenecks can both affect and be affected by surrounding traffic- and accounting for their interactions accurately and efficiently from a computational perspective. In this section, we describe how a fast algorithm based on the Lax-Hopf algorithm outlined earlier can be used to compute the solution associated with multiple moving (and fixed) bottlenecks.

In order to derive multiple internal conditions associated with several (active) moving bottlenecks, while accounting for their interactions and different properties, we propose a strategy based on two important properties of the solutions to Hamilton-Jacobi equations: the existence of a domain of influence for each internal condition, and the inf-morphism property of solutions.

*3.2. Inf-morphism and domains of influence*

The inf-morphism property implies that we can dynamically update the number of moving bottlenecks considered in the simulation problem, without having to re-compute the solution entirely. The domains on which the solution has to be re-computed are the domains of influence of the bottlenecks, which are the set of points that can be reached by characteristics with speeds ranging from $-w$ to $v$, and originating on the internal condition (moving bottleneck trajectory).

Indeed, for each position of the moving bottleneck, it is possible to identify its region of influence in the space and time dimension delimited by the congested and free-flow speed in the triangular fundamental diagram (respectively equal to $w$ and $v$). Whenever the moving bottleneck $i$ at the position $(x_i, y_i)$ enters in the domain of influence of the moving bottleneck j at the position $(x_j, y_j)$, the derivation of internal conditions has to be performed along the trajectory of moving bottleneck $j$. This stepwise computation can be repeated back and forth among several moving bottlenecks until the simulation is completed. Algorithm 2 summarizes this process, for an arbitrary number $n_b$ of moving bottlenecks with (possibly distinct) maximum speeds $v_{mb,i}$, entering the road at $(x_{1,i}, t_{1,i})$ and leaving at $x_{2,i}$.

*3.3. Passing bottlenecks and Zeno effect*

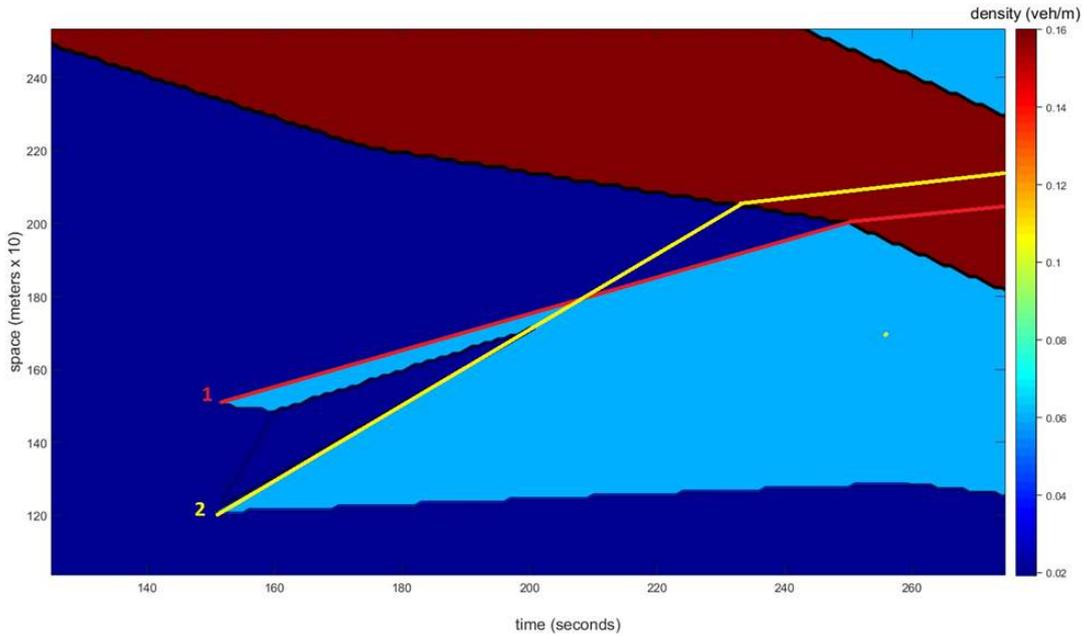

Figure 3: space-time-density diagram showing the simulation of two overtaking moving bottlenecks

The model used in this study corresponds to the coupling of a PDE (the LWR equation) with an ODE describing the evolution of the slow vehicles. The latter is hybrid, in the sense that the slow vehicles can only be in the three possible states outlined earlier. As with all hybrid systems, the dynamics can sometime exhibit the Zeno effect (Johansson et al., 1999). An execution of a hybrid system is called Zeno, if it takes infinitely many discrete transitions (and therefore computational loops) to solve the problem over a finite time horizon. In the present situation, the Zeno effect arises when bottlenecks are passing each other, as illustrated in Figure 3. In this situation, their respective domains of influence impose an upper bound on the time step used to update the position of each bottleneck (to ensure that the final position of one bottleneck is always outside the domain of influence of the other), and this upper bound becomes infinitely small as their paths come to intersect. This effect complicates the implementation of the algorithm as it can lead to infinite loops, if we want an exact solution. To solve this problem, we adopt a constant time step for the computation of the trajectories associated with moving bottlenecks: this allows the execution to be complete over a finite (and bounded above) number of steps, at the cost of computational accuracy, since this introduces an approximation of the behavior of both bottlenecks when they intersect each other.

### 3.4. Algorithm

The corresponding pseudo-code is shown as Algorithm 2 below.

Algorithm 2: Pseudo-code implementation for the computation of internal conditions associated with $n_b$ multiple moving bottlenecks

| | |
|---|---|
| **Input: $(x_{1,i}, t_{1,i}, v_{mb,i}) \ \forall \ i \in [1, n_b]$** | *Input initial position and time, final position, and performance characteristics of moving bottlenecks* |
| **Input: $(v, k_c, k_j, n_l)$** | *Input list of fundamental diagram parameters* |
| **Input T** | *Input simulation time horizon* |
| $q_{r,i} = (v - v_{mb,i}) * k_c * (n_l - 1)/n_l;$ | *Derive maximum passing rate for each moving bottleneck* |
| $b = \{1, \ldots, n\_b\}$ | *Initialize bottlenecks list* |
| $t_i = t_{1,i}$ | *Current time for each bottleneck* |
| $x_i = x_{1,i}$ | *Current position for each bottleneck* |
| while $b \neq \emptyset$ | *While some bottleneck vehicles are still on the computational domain* |
| | *Eliminate bottlenecks that have left the computational domain* |
| for $i \in b$ <br>   if $t_i > T$ or $x_i \geq x_{max}$ <br>   $b = b \setminus \{i\}$ <br>   end <br> end | *Initialize list of bottlenecks that are not influenced by others (and that can thus be computed)* |
| $l = b$ | *Compute list list of bottlenecks that are not influenced by others* |
| for $i \in b$ <br>  for $j \in b \setminus \{i\}$ <br>   if $(t_i, x_i) \in D_j$ ($D_j$ domain of influence of bottleneck j) <br>    $l = l \setminus \{i\}$ <br>   end <br>  end <br> end | |
| | *If there are moving bottlenecks that are not influenced by others* |
| if $l \neq \emptyset$ | *Propagate these bottlenecks according to Algorithm 1* |
|   while $l \neq \emptyset$ <br>     pick $i \in l$ <br>     **compute propagation of moving bottleneck** *(according to Algorithm 1)* <br>     $l = l \setminus \{i\}$ <br>   end while <br> else | *If none of the bottleneck is outside of the zone of influence of all others* (bottleneck intersection) |
| **Identify intersecting bottlenecks, and propagate them approximately according to Algorithm 1** | |
| end | |
| end while | |

An example of the application of Algorithm 2 is illustrated in the sequence of simulations in Figure 4, where the different steps of the computation of the solution are shown for a pair of moving bottlenecks (indicated by letters A and B). First, the trajectory of moving bottleneck A and its impact on traffic are computed till it intersects the domain of influence of moving bottleneck B (Figure 3a), indicated by the green lines (Phase 1). Likewise, the solution associated with moving bottleneck B can be computed till it reaches the domain of influence of moving bottleneck A at its last position in the space and time (Figure 3b). The procedure is repeated back and forth till the moving bottleneck leaves the road (Figure 3c) or the simulation ends (Figure 3d).

A more general example involving ten moving bottlenecks having different speeds (Table 3) and three fixed bottlenecks (representing constant red cycles of a traffic light) on the same link used in the previous cases, is illustrated in Figure 5. The computational performance, consisting of 0.15 seconds for the computation of the internal conditions and 2.2 seconds for the computation of the solutions on a rectangular grid of resolution one second and ten meters, confirms that the algorithm can handle complex scenarios very efficiently.

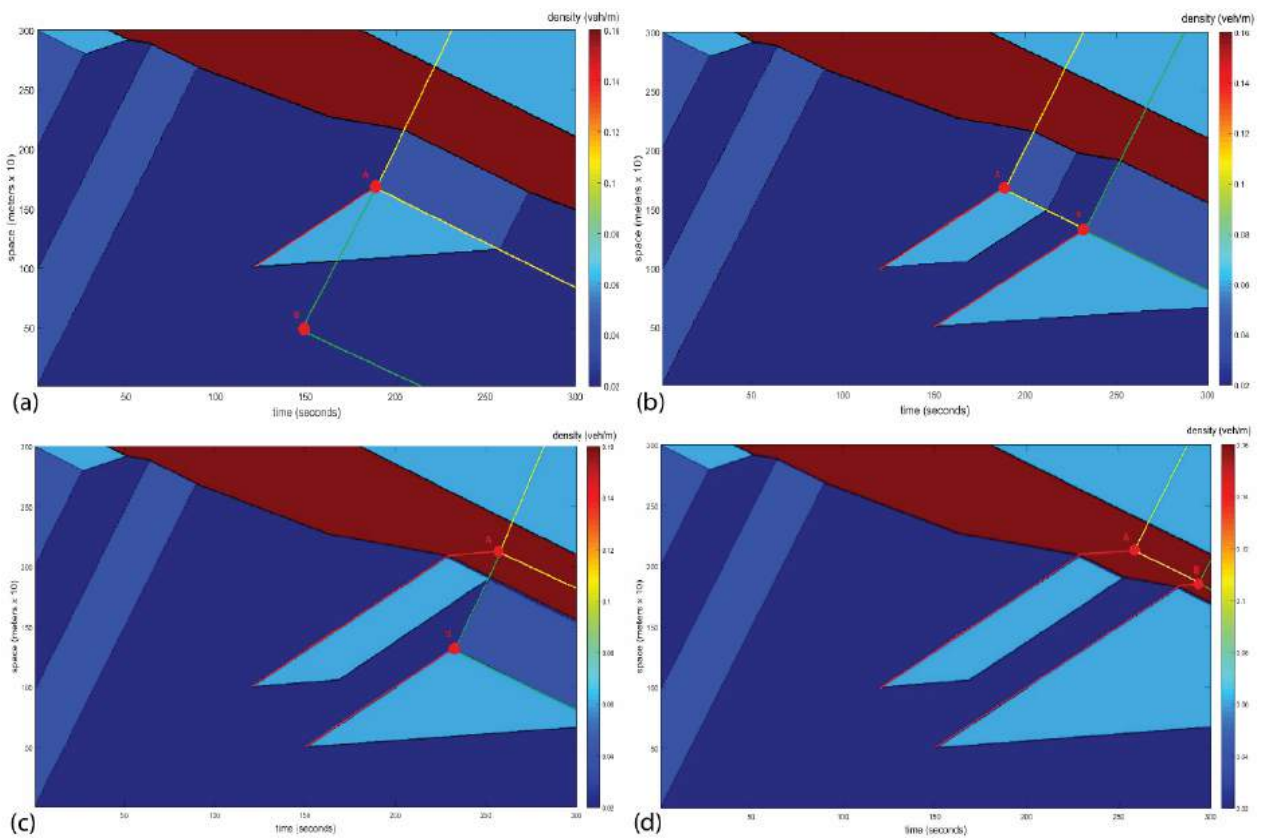

Figure 4: Example propagation of two moving bottlenecks. In (a), moving bottleneck A intersects the domain of moving bottleneck B. The algorithm continues propagating both bottlenecks (b, c) until both bottlenecks have left the computational domain (d)

Table 3: Moving bottleneck features

| $i$ | 1 | 2 | 3 | 4 | 5 | 6 | 7 | 8 | 9 | 10 |
|---|---|---|---|---|---|---|---|---|---|---|
| $x_1$ | 2000 | 1000 | 1000 | 1600 | 1200 | 2000 | 800 | 1500 | 1500 | 1000 |
| $t_2$ | 60 | 20 | 50 | 150 | 120 | 220 | 180 | 270 | 330 | 320 |
| $v_{max}$ | 5 | 8 | 10 | 10 | 8 | 12 | 10 | 8 | 5 | 5 |

## 4. Application to schedule and traffic signal optimization

In order to show the capabilities of the algorithm, two different optimization problems involving several moving and fixed bottlenecks are investigated. The first one consists in the maximization of the total outflow on a two-lane road where the decision variables include adjustments of the entry times of eight buses, and the duration of signals' cycles and green phases. In the second one, on the same typology of scenario, the same decision variables are varied to minimize buses' delay.

In these tests, we assume that the initial and boundary conditions are known (and can be arbitrary). Though real-world applications of signal optimization involve multiple connecting links, we choose a single link (with two traffic signals) for this test. The extension of the present application to networks is straightforward, and involves the use of dynamic boundary conditions for each link, these conditions being linked together through junction models.

### 4.1. Formulation and computation approach

The approach adopted to solve both optimization problems is based on the search heuristic technique of the genetic algorithms (GA), which are an optimization strategy where a set of randomly generated solutions (initial population) is improved by means of an iterative procedure. This iterative process consists of selecting the best performing solutions (parents) and "breeding" them to create new "generations" of solutions (children), until an optimal solution is found. During the breeding process mutations are randomly applied (i.e. random changes in a solution that occur at a predetermined probability) in order to maintain a higher diversity in the population. The interested reader may refer to Yang (2010) for a more detailed explanation of this heuristic.

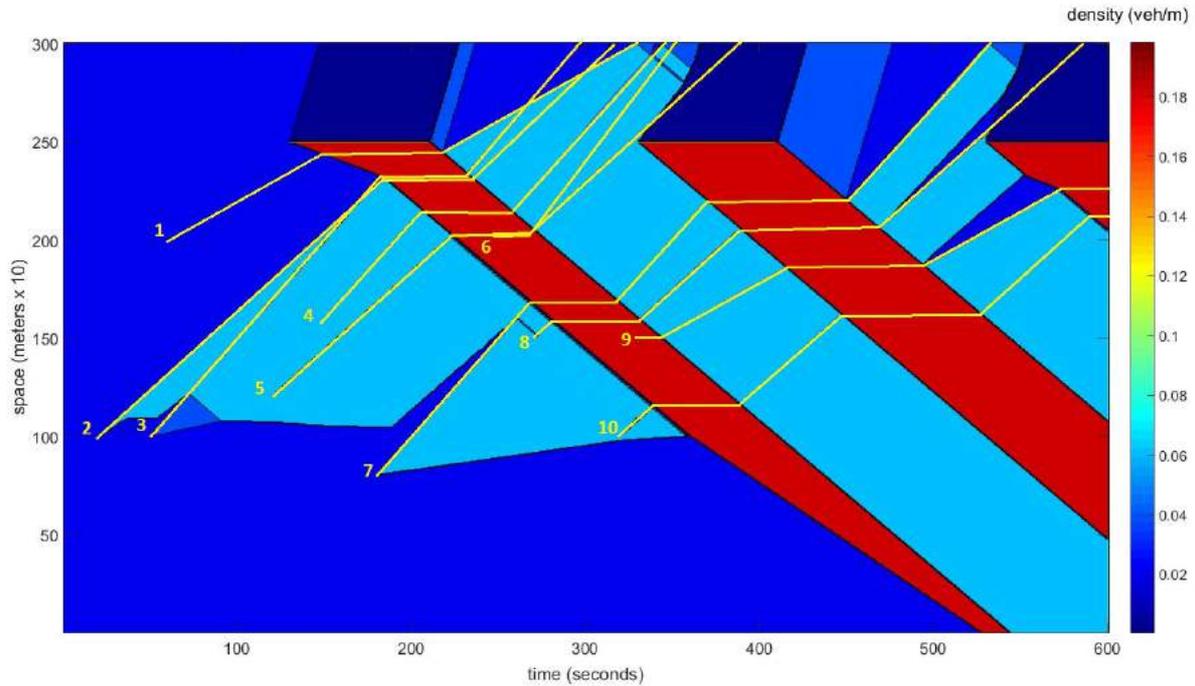

Figure 5: Example of simulation of several fixed and moving bottlenecks having different maximum speeds

As Teklu et al. (2007) observe, GA are very flexible as they do not require any knowledge of the gradient and they can avoid getting stuck in local optima. Furthermore, GA are more suitable to simulation-based frameworks, as derivative-based optimization methods require the knowledge of the analytical form of the problem. For these reasons GA have been already employed to solve different joint traffic control and assignment problems (Foy et al., 1992; Lee and Machemehl, 1998; Yin, 2000; Ceylan and Bell, 2004; Teklu et al., 2007).

In both the optimization problems, the objective functions depend on the entry times of i moving bottlenecks $t_i$ and the signal setting variables $\varphi_j = (c, \phi)$, where $c$ and $\phi$ correspond respectively to the cycle time length and the green time length of traffic signal j.

In the first problem, we optimize the total outflow, that is, the cumulative number of vehicles at the downstream end of the link, $N_d$:

$$\max_{t_i, \varphi_j} N_d(t_i, \varphi_j)$$

subject to:
1. $t_{i,min} \leq t_i \leq t_{i,max}$  $\forall i \in I$            moving bottlenecks' entry times constraints
2. $\varphi_j(c, \phi) \begin{cases} c_{min} \leq c \leq c_{max} \\ \phi_{min} \leq \phi \leq \phi_{max} \end{cases}$  $\forall j \in J$        cycle times and green times constraints

In the second optimization problem, which is characterized by the same decision variables and constraints of the first optimization problem (1,2), the total delay $D$ of moving bottlenecks is minimized. This yields to the following formulation:

$$\min_{t_i, \varphi_j} D(t_i, \varphi_j)$$

Where the total delay is given by the sum of each moving bottleneck's delay, calculated as difference between a theoretical arrival time if the moving bottleneck travels at its maximum speed ($v_{mb,i}$) and its actual arrival time ($t_{real,i}$) derived from the simulation:

$$D(t_i, \varphi_j) = \sum_{i \in I} \frac{(x_{2,i} - x_{1,i})}{v_{mb,i}} - t_{real,i}$$

with $x_{1,i}$ and $x_{2,i}$ corresponding to the input entry and exit points of each moving bottleneck i.

For the computation of GA, each solution is identified as a vector of values corresponding to each decision variables. The main algorithmic steps of the GA consist of:

1. Generate the initial random population of solutions corresponding to different moving bottlenecks' entry times and signal timings (satisfying the problem constraints)
2. Calculate the fitness functions (according to the optimization problem) by means of simulation for each solution
3. Perform breeding process by means of tournament selection to select the parents and one-point crossover to generate the children
4. Perform random mutation of children with probability Pm, by randomly replacing some their decision variables' values (satisfying the problem constraints)
5. If the maximal generation number is reached, the solution with the highest fitness is adopted as optimal solution of the problem. Else, return to the second step.

*4.2. Numerical results*

In both optimizations problems the initial scenario consists of a two-lanes link of 3000 meters length, with two traffic lights at the positions $y_1 = 500\ m$ and $y_2 = 2000\ m$. Both traffic lights are characterized by an initial cycle time of 200s and green time of 120s. The boundary and initial conditions are reported respectively in Table 1 and Table 2. Eight buses (moving bottlenecks) enter and leave the link at different points (Table 3), and they are characterized by the same maximum speed $v_{mb} = 20 \ m/s$. The following parameters characterizing the triangular fundamental diagram are chosen: $u = 30 \ m/s$, $k_c = 0.02 \ veh/m$, $k_j = 0.1 \ veh/m$. The results of the simulation are illustrated in Figure 6.

Table 4: boundary conditions

| $i$ | $t_{i-1}$ | $t_i$ | $q_{up}^i$ | $q_{down}^i$ |
|---|---|---|---|---|
| 1 | 0 | 40 | 1.2 | 0.5 |
| 2 | 40 | 300 | 1.2 | 1.0 |

Table 5: initial conditions

| $i$ | $x_{i-1}$ | $x_i$ | $k_{ini}$ |
|---|---|---|---|
| 1 | 0 | 3000 | 0.04 |

Table 6: moving bottleneck conditions

| $i$ | 1 | 2 | 3 | 4 | 5 | 6 | 7 | 8 |
|---|---|---|---|---|---|---|---|---|
| $x_1$ | 1000 | 1000 | 2000 | 1000 | 2000 | 1000 | 2000 | 2000 |
| $t_1$ | 20 | 50 | 60 | 120 | 150 | 180 | 220 | 270 |
| $x_2$ | 3000 | 3000 | 3000 | 2000 | 3000 | 2000 | 3000 | 3000 |

As to the constraints of the optimization problem, the entry times of buses can be moved ahead or postponed up to 10 seconds, and cycle lengths and green times can be shrunk or extended by a factor of 0.1. In a network optimization problem, the boundary flows would become additional decision variables, and could similarly be constrained to be in some interval.

The evolution of the fitness function through several generations in both optimization problems suggests the designed GA converges to a maximum/minimum rather quickly (between the 10[th] and 15[th] iteration) regardless of the input parameters used (see Figure 7). From a computational perspective this result implies that the optimization algorithm can be limited to fewer generations without compromising the quality of the results. More interestingly, as several tests concerning the size of the population of the GA show, the quality of the solution does not seem to be largely affected in both optimization problems (Table 4 and Table 5). Indeed, in the case of maximization of outflow, increasing the size of the population from 30 to 100 seems to improve only by 3 percent the fitness function (on average), at the expense of the computational time, which is more than tripled. Similarly, GA with a population of 100 solutions would yield to an average improvement of 10 percent compared to GA with a population of 30 solutions, but with a computation times more than three times higher.

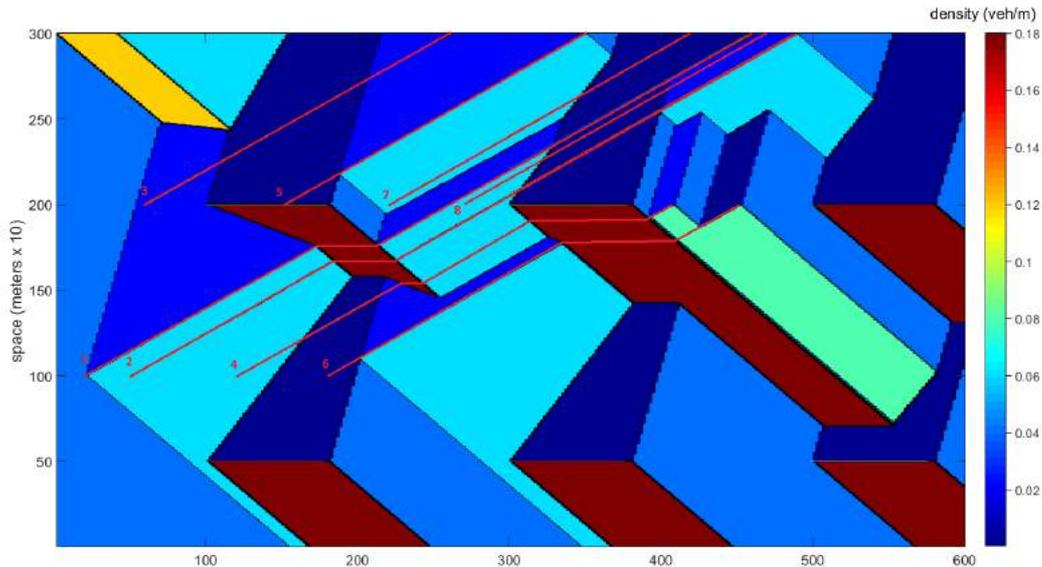

Figure 6: : Space-Time-Density diagram for the original scenario before optimization

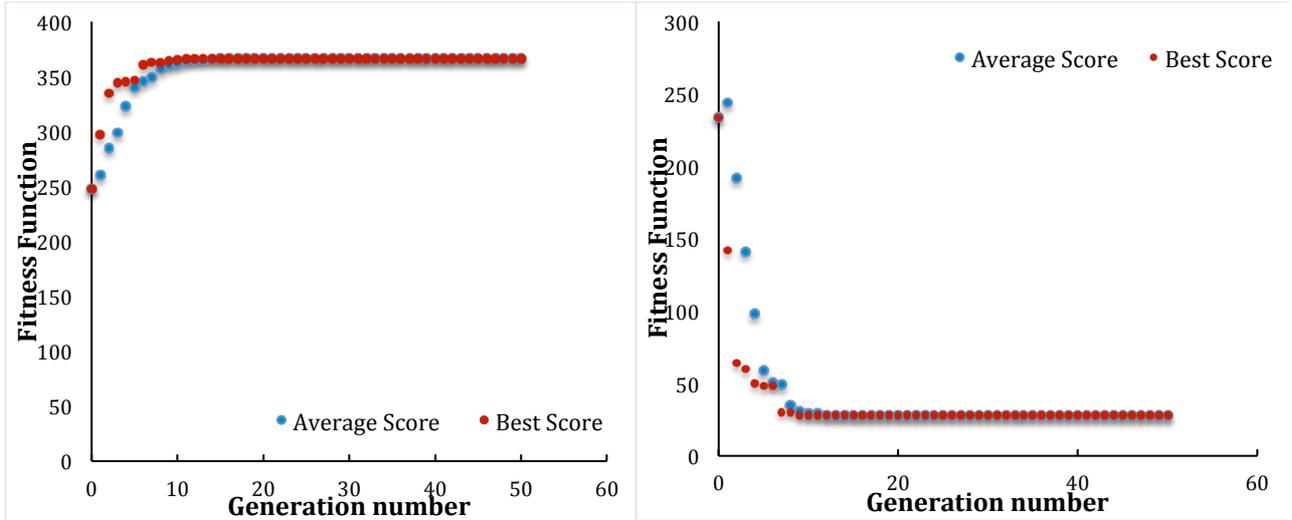

Figure 7: Evolution of the fitness function over several generations in the two optimization problems[†]

The optimal solutions corresponding to both optimization problems are displayed in Figure 8. As shown in Figure 8a, the outflow can be increased, not only by varying the cycle lengths of the two traffic lights, but also by reducing the delays caused by the moving bottlenecks (for example bus 1 and bus 5). When we minimize the bus delay instead (Figure 8b), mainly the buses' trajectories and the settings of the more downstream traffic light are optimized to reduce buses' wait at the traffic lights.

The required computation time to calculate the objective function for each of the solutions created during the breeding procedure and accomplished by means of traffic simulation, varies between 0.05 and 0.07 seconds.

Overall, this analysis suggests that, for this optimization problem involving eight moving bottlenecks and two traffic lights, it is possible to obtain fairly satisfactory solutions in a few tens of seconds, thanks to the fast algorithm introduced earlier, which we use to compute the solutions associated with multiple moving bottlenecks.

Table 6:

|  | Best score | Solution improvement (%) | Computation time |
| --- | --- | --- | --- |
| Population 30 | 335 | 26 | 41.1 |
| Population 50 | 344.3 | 28 | 67.4 |
| Population 100 | 346.2 | 28.4 | 133.9 |

Table 7:

|  | Best score | Solution improvement (%) | Computation time |
| --- | --- | --- | --- |
| Population 30 | 70.6 | 70 | 39.6 |
| Population 50 | 48.3 | 79 | 70.6 |
| Population 100 | 46.8 | 80 | 140.5 |

---

[†] The plotted evolutions correspond to the best results obtained for a simulation with population of 100 solutions, tournament size of 10 and mutation rate of 0.1

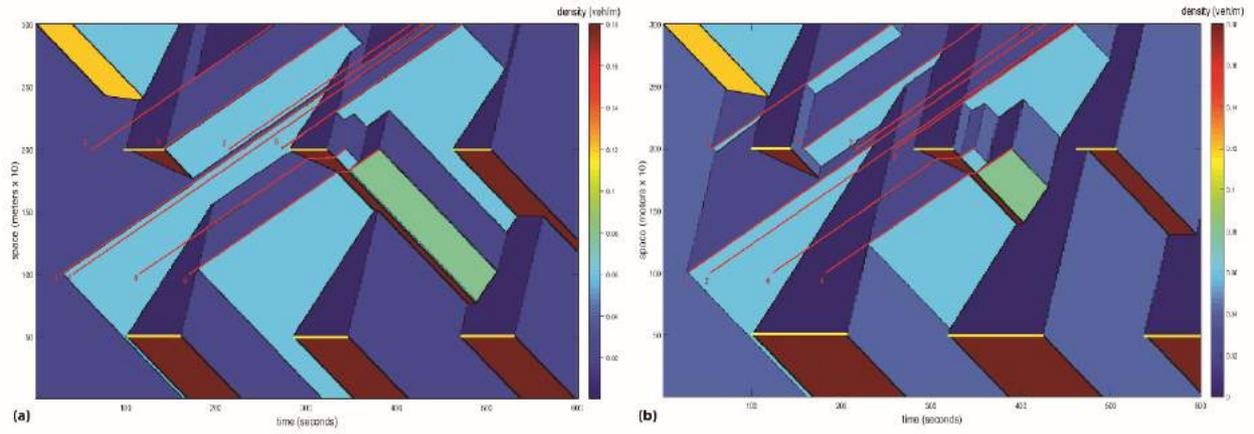

Figure 8: Space-time-density diagram corresponding to the results of the first optimization problem (a) and the second optimization problem (b). Buses' trajectories are shown as red lines, and red traffic lights are indicated by yellow lines.

## 5. Conclusions and future work

In this study, we proposed a new semi-analytic numerical scheme that can be used to compute the solutions within the LWR traffic flow model given initial, upstream and downstream boundary conditions, and an arbitrary number of moving bottlenecks, which can be associated with different types of vehicles. The main feature of this scheme is the ability to endogenously and efficiently account for the interaction between several moving (and fixed) bottlenecks and surrounding traffic, without relying on continuous flow or multi-class approximations.

This numerical scheme is based on a Hamilton-Jacobi formulation of the LWR model, and results from the properties of the solutions to Hamilton-Jacobi equations, and in particular the inf-morphism property. Being semi-analytic, it is very accurate (though not exact due to the piecewise linear approximation of the trajectories of the moving bottlenecks), and very fast, since it allows one to determine the trajectories of all moving bottlenecks without having to compute the solution on the entire computational domain, making it very adapted to optimization problems.

We demonstrated the benefits of this algorithm by taking advantage of its high computational performance in two different optimization problems where we simultaneously optimized the schedule of eight buses and the timing of two traffic signals on a single stretch of road. These problems could be solved using known heuristics in a few dozens of seconds with a regular computer.

Future work will deal with the extension of this framework to the optimization of schedules and traffic signals simultaneously over road networks. In the network case, the boundary conditions on each domain are not fixed, and the moving bottlenecks can propagate between domains. Extending the proposed framework to networks will require the estimation of boundary flows of each of the links by means of junction models, which could slightly increase the necessary computational time to compute the solution.